\documentclass[a4paper, 11pt, DIV=12,abstract=on]{scrartcl}
\usepackage[T1]{fontenc}
\usepackage[utf8]{inputenc}
\usepackage[english]{babel}
\usepackage{amssymb}
\usepackage{amsmath}
\usepackage{amsfonts}
\usepackage{mathdots}
\usepackage{graphicx}
\usepackage[inline]{enumitem}
\usepackage{url}

\newtheorem{definition}{Definition}[section]
\newtheorem{theorem}{Theorem}[section]
\newtheorem{lemma}{Lemma}[section]
\newtheorem{corollary}{Corollary}[section]
\newcommand{\Proof}[0]{\smallskip\noindent\normalfont\underbar{Proof}: }
\newcommand{\Black}[0]{$~\blacksquare$}

\newcommand{\Col}{\overleftarrow{C}}  
\newcommand{\Cor}{\overrightarrow{C}} 
\newcommand{\Cgb}{{\boldsymbol{\mathrm{C}}}_b} 
\newcommand{\Cgh}{{\boldsymbol{\mathrm{C}}}_h} 
\newcommand{\Cgr}{{\boldsymbol{\mathrm{C}}}_{\overrightarrow{C}}} 
\newcommand{\Fob}{\boldsymbol{F}_b}
\newcommand{\Foh}{\boldsymbol{F}_h}
\newcommand{\Atr}{{\boldsymbol{A}}_{\overrightarrow{C}}} 

\newcommand{\setE}{\mathbb{E}} 
\newcommand{\setN}{\mathbb{N}} 
\newcommand{\setO}{\mathbb{O}} 
\newcommand{\setY}{\mathbb{Y}} 
\newcommand{\setZ}{\mathbb{Z}^+} 

\begin{document}
\author{Heinz Ebert}
\title{A Graph Theoretical Approach to the Collatz Problem}
\maketitle
\begin{abstract}
\noindent Andrei et al. have shown in 2000 that the graph $\boldsymbol{\mathrm{C}}$ of the Collatz function starting with root $8$ after the initial loop is an
infinite binary tree $\boldsymbol{A}(8)$. According to their result they gave a reformulated version of the Collatz conjecture: the vertex set
$V(\boldsymbol{A}(8))=\setZ$.

In this paper an inverse Collatz function $\Cor$ with eliminated initial loop is used as generating function of a Collatz graph $\Cgr$. This graph can be
considered as the union of one forest that stems from sequences of powers of 2 with odd start values and a second forest that is based on branch 
values $y=6k+4$ where two Collatz sequences meet. A proof that the graph $\Cgr(1)$ is an infinite binary tree $\Atr(1)$ with vertex set $V(\Atr(1))=\setZ$ completes the paper.


\smallskip \noindent \textbf{MSC-Class:} 11B83, 05C05, 05C63
\end{abstract}

\section{The Collatz function and conjecture}

Let $\setN$ be the set of nonnegative integers and $\setZ$~be the positive integers, then the Collatz problem relates to the Collatz map $\Col$:
$\setZ\to\setZ$:
\begin{equation}
	\text{$\Col (n)=$}
	\begin{cases}
		n/2 & \text{if $n\equiv 0\ (mod\ 2)$,$\quad \Col(n)\in \setZ$} \\
		3n+1 & \text{if $n\equiv 1\ (mod\ 2)$,$\quad \Col(n)\equiv 4\ (mod\ 6)$.}
	\end{cases}
\end{equation}
\noindent The famous 3n+1 or Collatz conjecture now states that for any $n\in \setZ$ there exists a $k\in \setN$ such that:
\begin{equation*}
	\text{$\Col^{(k)}(n)=1,~~ \bigl\lbrack\Col^{(0)}(n)=n$ and $\Col^{k}(n)=\Col\circ\Col^{k-1}(n)\bigr\rbrack$.}
\end{equation*}

\noindent The conjecture excludes the existence of other loops than the trivial terminal cycle $(1, 2, 4, 1,\ldots)$ and of any divergent sequences.

\section{The Collatz tree and a modified conjecture} 

Most papers deal with the dynamics of the Collatz function $\Col$ or modified versions of it while pure graph theoretical aspects have seldom been considered.
Some exceptions are Andaloro {\footnotesize \fbox{1}}, Andrei et al. {\footnotesize \fbox{2,3}}, Laarhoven and de Weger {\footnotesize \fbox{6}}, Lang
{\footnotesize \fbox{7}} and Wirsching {\footnotesize \fbox{8}}.

\noindent Andrei et al. {\footnotesize \fbox{2,3}} examined a graph $\boldsymbol{\mathrm{C}}$ of the Collatz function and showed that a subgraph of
$\boldsymbol{\mathrm{C}}$ with the vertex set $V\subseteq \setZ-\{1,~2,~4\}$ and the value 8 as root is an infinite binary tree $\boldsymbol{A}(8)$. Therefore they called it Collatz tree. According to this result they reformulated the Collatz conjecture to be:

\medskip \ \ \ The vertex set of the Collatz tree $\boldsymbol{A}(8)$ is $V=\setZ-\{1,2,4\}$.

\medskip \noindent Their conclusions also lead to the fact that every $n>4$ could be the root of a Collatz tree $\boldsymbol{A}(n)$. Then they concentrate on
infinite chain subtrees which are characterized by values which are divisible by $3$. Graphs without these chain subtrees are called pruned Collatz graphs
{\footnotesize \fbox{8}}. This approach leads to infinite sets of start numbers whose sequences converge at 1.

\section{The inverse Collatz function} 

\noindent Let the set $\setY=\{n>4|n\equiv 4\ (mod\ 6)\}\subset \setZ$, then the inverse Collatz map $\Cor$: $\setZ \to \setZ$ is:
\begin{equation}
	\text{$\Cor (n)=$}
	\begin{cases}
		2n & \text {if $n\in \setZ$,\ \ \ $\Cor(n)\equiv 0\ (mod\ 2)$} \\
		(n-1)/3 & \text {if $n\in \setY$,\ \ \ \ \ $\Cor(n)\equiv 1\ (mod\ 2)$.}
	\end{cases}
\end{equation}
Although the two operations of the Collatz function $\Col$ have the above unique inverses in the definition of $\Cor$, the function $\Cor$ itself is not unique.
This is because $\setY$ is a proper subset of $\setZ$. This leads to the fact that every $y\in \setY$ always has two descendants. It is obvious that the
operation $2n$ simply continues its current sequence while the operation $(n-1)/3$ results in an odd number and starts a complete new sequence. 
Therefore we call the numbers $y$ \emph{branch values}. As $4$ is such a branch value we excluded $4$ from the set $\setY$ to avoid the otherwise inevitable initial loop
$(\boldsymbol{1},~2,~4,~\genfrac{}{}{0pt}{}{\boldsymbol{1},\ldots }{8,\ldots })$. 

\section{The Collatz graph of the inverse Collatz function} 

In 1977 Lothar Collatz remarks in a paper on the use of graph representations to study iteration problems of functions ${f}$: $\setZ\to\setZ$ but he did not
consider the 3n+1 problem therein {\footnotesize \fbox{4}}. His idea was to picture such dynamical systems by infinite graphs of the following kind:
\begin{definition}
	\normalfont Let $\bigl(n, f(n)\bigr)\in\setZ$, then an infinite Collatz graph is generally defined by:
	\begin{equation}
		\text{${\boldsymbol{\mathrm{C}}}_f(V_f,E_f)=$}
		\begin{cases}
			V_f= \setZ & \text{\normalfont{the set of vertices}} \\
			E_f= \bigl\langle n,f(n)\bigr\rangle~~\bigl(n,f(n)\bigr)\in V_f & \text{\normalfont{the set of \underbar{directed} edges.}}
		\end{cases}
	\end{equation}
\end{definition}

\noindent There are important differences between \emph{normal} graphs and Collatz graphs:
\begin{enumerate}[nosep]
	\item The vertex set and the results of the generating function $f$ are restricted to the set $\setZ$.
	\item Vertices and their labels are indistinguishable.
	\item The map $f$ determines the set of edges and their direction $n\to{f(n)}$.
	\item The map $f$ enforce the properties of the vertices/labels.
\end{enumerate}

\medskip \noindent An example for the above point 4 are the numbers $y\in\setY$. The property induced by the maps $\Cor$ and $\Col$ is that all $y\equiv 4\ (mod\ 6)$. But $\Col$ named now as Collatz backward function ignores that these numbers are \emph{branch values}. The inverse Collatz function $\Cor$ is therefore much more appropriate as a generating function of a graph and so we use from now on the map $\Cor$ named as Collatz forward function for the construction of our Collatz graphs exclusively. Since $\Cor$ is the inverse map of $\Col$ we even relax the demand of Definition 4.1 that the edges have to be directed. Thus we define the common graph for both Collatz functions as:

\begin{equation}
	\text{$\Cgr=$}
	\begin{cases}
		V_{\Cor}=\setZ & \text{the set of vertices} \\
		E_{\Cor}=\bigl(n,\Cor(n)\bigr) \quad \bigl(n,\Cor(n)\bigr)\in V_{\Cor} & \text{the set of \underbar{undirected} edges.}
	\end{cases}
\end{equation}
\noindent and the Collatz conjecture reads now:

\smallskip \ \ \ The graph $\Cgr(1)$ is an undirected infinite binary tree $\Atr(1)$ with the vertex set $V=\setZ$.

\section{The Collatz graph $\protect \Cgr$ as union of two infinite forests} 
\noindent Diestel defines a forest as: \emph{A graph without circles is called a forest. A connected forest is a tree. Thus a forest is a graph whose components are trees} {\footnotesize\fbox{5}}.

\subsection{The forest $\protect \Foh$}

\noindent We now show what happens if we repeatedly apply the operation $n'=2n$ of $\Cor$ to all odd start numbers $o\in\setO=\{n>0|n\equiv 1~(mod~2)\}$. The inverse operation is $n'=n/2$ of $\Col$ applied to any even number $n\in\setE=\{n>0|n\equiv 0~(mod~2)\}$ until $n'$ is odd.

\begin{figure}[htbp]
	\centering
	\includegraphics[scale=0.4]{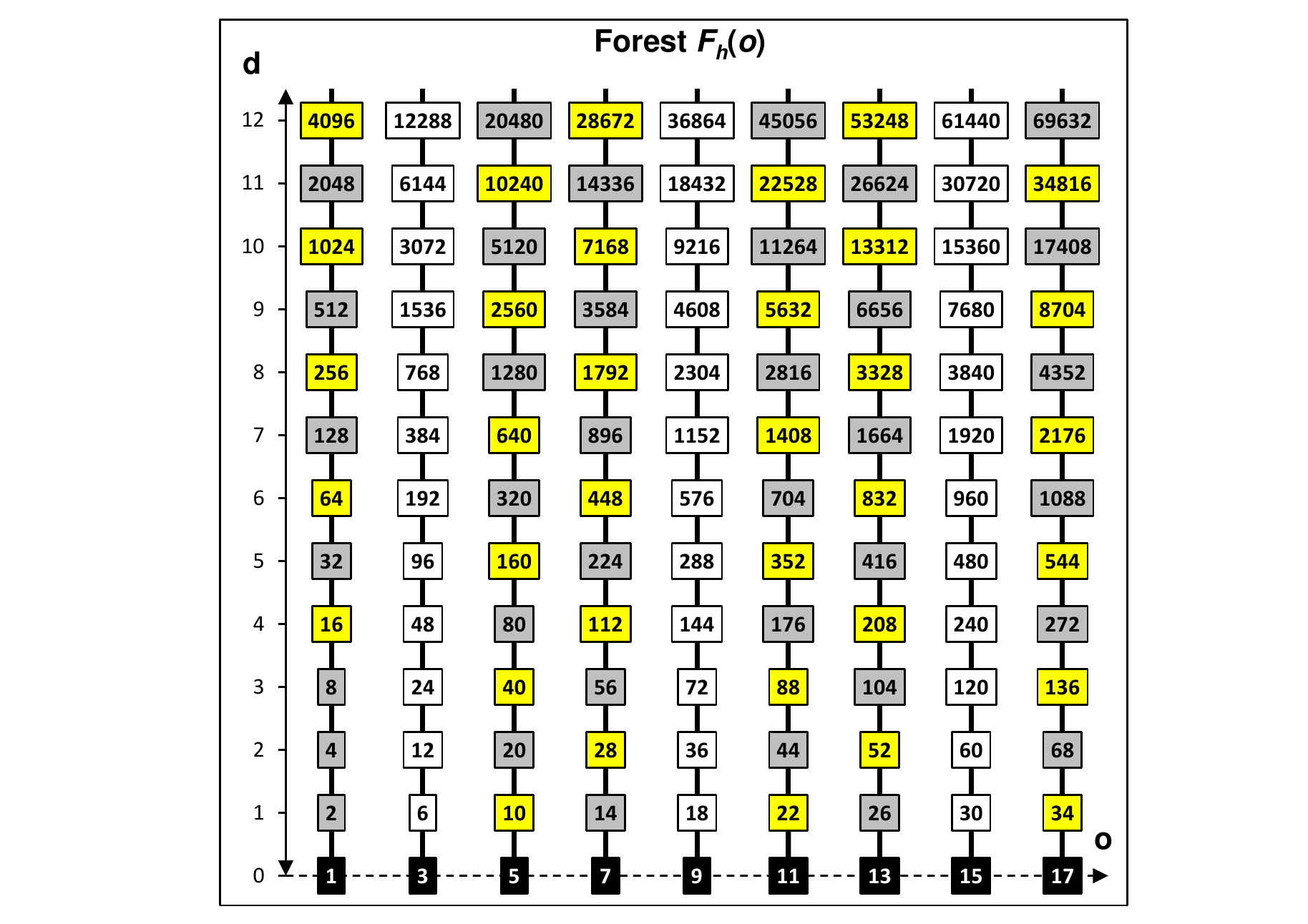}
	\caption{Grid graph of the Forest $\protect \Foh$. The generating function $h(o,d)=o\cdot 2^d$ dictates the colors 
		indicating the properties of the nodes: $v\equiv 1~(mod~2)$ black, $v\equiv 4~(mod~6)$ yellow, $v\equiv 2~(mod~6)$ grey, 
		$v\equiv 0~(mod~6)$ white.}
	\label{fig:Fh}
\end{figure}

\begin{theorem}
	\noindent Let $o\in\setO$ and $d\in\setN$, then with $o\to\infty$ and $d\to\infty$ the Collatz graph $\Cgh$ 
		generated by the function $h(o,d)=o\cdot 2^d$ is an infinite forest $\Foh$ 
		of distinct infinite trees ${\boldsymbol{A}}_h(o)$ with the set of vertices $V(\Foh)=\setZ$.

	\Proof For any fixed $O\in\setO$ and $d\to\infty$ the infinite sequence $h(O,d)=O\cdot 2^d$ resembles 
		a single infinite tree ${\boldsymbol{A}}_h(O)$ without any branches. 
		Thus with $o\to\infty $ we get a set of unconnected infinite trees: the forest $\Foh$ with the set of edges
		$E(\Foh)=\bigl\{e|e=o\cdot 2^d,~o\cdot 2^{d+1}\bigr\}$ (Figure~\ref{fig:Fh}).\\
		For d=0 the codomain of $h(o,0)$ is the set $\setO$ and for $d>0$ the codomain of $h(o,d)$ is the set $\setE$. 
		The set of vertices of $\Foh$ is $V(\Foh)=\setO\cup\setE=\setZ$. \Black
\end{theorem}

\begin{corollary}
Obviously all vertices $o\in\setO$ as roots of the trees ${\boldsymbol{A}}_h(o)$ have one incident edge and all nodes $v\in \mathbb{E}$ have two incident edges.\end{corollary}

\subsection{The forest $\protect \Fob$}

\begin{figure}[htbp]
	\centering
	\includegraphics[scale=0.4]{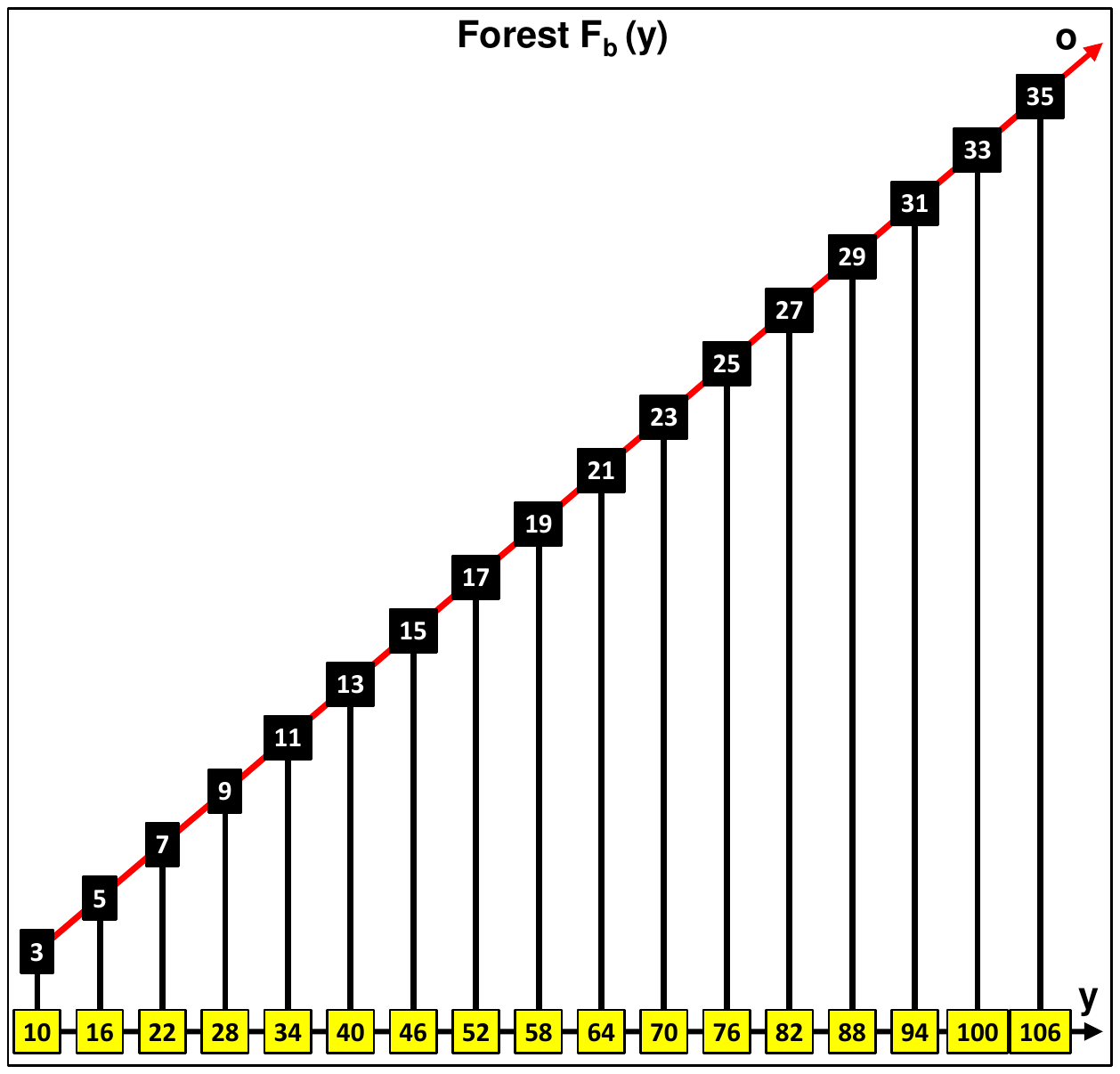}
	\caption{Grid graph of the Forest $\protect \Fob$. The generating function is $b(y)=(y-1)/3$ and the properties of the vertices are: 
		$v\equiv 1\ (mod\ 2)$  black, $v\equiv 4\ (mod\ 6)$ yellow.}
	\label{fig:Fb}
\end{figure}

\noindent Now we exclusively apply the operation $o=(y-1)/3$ of $\Cor$ to all branch numbers $y>4$. The inversion is the operation $y=3o+1$ of $\Col$ applied to all numbers $o>1$.
\begin{theorem}
	Let $o\in\setO$, $y\in\setY$ and the map $b$: $\setY\to\setO$: $b(y)=(y-1)/3$, then with $y\to\infty$ 
	the Collatz graph $\Cgb$ is an infinite forest $\Fob$ of distinct trees ${\boldsymbol{A}}_b(y)$.

	\Proof $E(\Cgb)=\{e|e=(y, o)\}$ and $V(\Cgb)=(\setY\cup\setO)\subset\setZ$. 
	Since all edges $e \in E(\Cgb)$ are different each edge $e$ represents a single tree ${\boldsymbol{A}}_b(y)$. 
	With $y\to\infty$ we get the forest $\Fob$ as set of infinitely many unconnected trees ${\boldsymbol{A}}_b(y)$ (Figure~\ref{fig:Fb}). \Black
\end{theorem}

\begin{corollary}
	Obviously all vertices $y\in\setY$ and $o\in\setO$ of the trees ${\boldsymbol{A}}_b(y)$ have one incident edge.
\end{corollary}

\subsection{Consequences of the union of \protect $\boldsymbol{\Foh}$ and \protect $\boldsymbol{\Fob}$} 

\begin{figure}[htbp]
	\centering
	\includegraphics[scale=0.4]{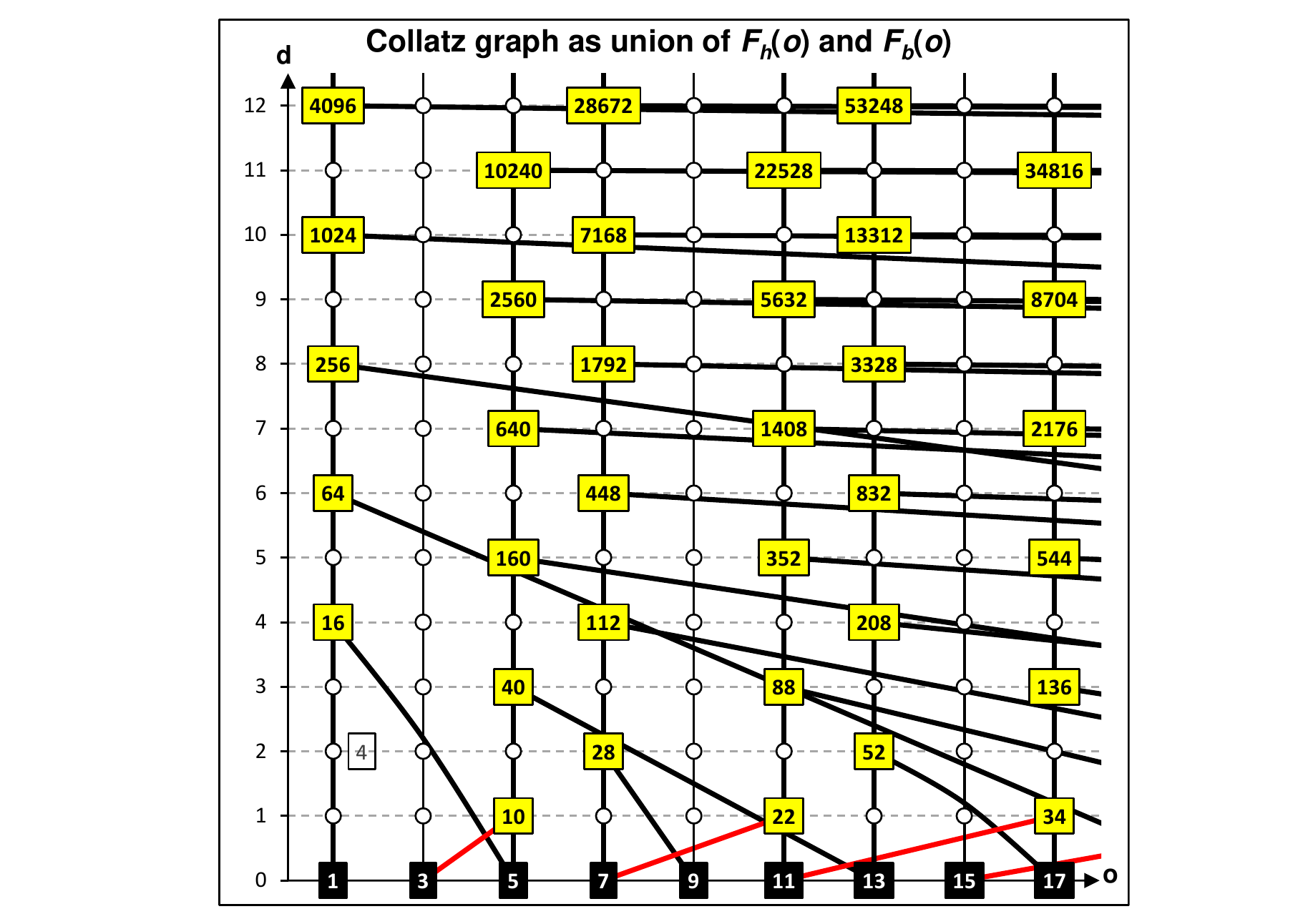}
	\caption{The grid graph $\protect \Cgr$. The forest $\protect \Foh$ rules the vertical and the forest $\protect \Fob$ the diagonal edges of this graph. 
		Circles represent nodes $v\equiv 2\ (mod\ 6)$, $v\equiv 0\ (mod\ 6)$ and $v=4$.}
	\label{fig:Ct}
\end{figure}

\noindent The separate application of operations of the generating functions $C$ and $\Cor$ split the Collatz graph $\Cgr$ into two different forests. The
re-union of $\Foh$ and $\Fob$ changes the sets of edges and the incidences of the nodes of both forests (Figure~\ref{fig:Ct}).
\begin{lemma}
	$E(\Foh)\cap E(\Fob)=\{0\}$.

	\Proof\\
	 Let $d\in\setN$, $o\in\setO$, $y\in\setY$, $E(\Foh)=\{e|e=({o\cdot 2}^d,~{o\cdot 2}^{d+1})\}$ and $E(\Fob)=\{e|e=(y,~o)\}$, 
	then all $e_{h,2}=o\cdot 2^{d+1}$ of $E(\Foh)$ are even and all $e_{b,2}=o$ of $E(\Fob)$ are odd and therefore all 
	edges of $E(\Foh)$ and $E(\Fob)$ are different. \Black
\end{lemma}

\begin{theorem}
	$\Cgr = \Foh\cup\Fob$.

	\Proof\\
	Because of Lemma 5.1 the union $E(\Foh)\cup E(\Fob)=E(\Cgr)$ introduces no multiple edges. 
	As $V(\Foh)=\setZ$ and $V(\Fob)\subset\setZ$ therefore $V(\Cgr)=V(\Foh)\cup V(\Fob)=\setZ$. \Black
\end{theorem}

\begin{theorem}
	All nodes $v\in V(\Cgr )$ have at most three incident edges.

	\Proof\\
	\noindent Due to Lemma 5.1 and Theorem 5.3 we can add and count the incident edges of $E(\Cgr)$:
	\begin{enumerate}[nosep]
		\item The root $v=1$ is no vertex of $\Fob$ and so only has one undirected edge $e=(1, 2)$.
		\item For all nodes $o>1$ there exist two undirected edges $(o,y)$, $(o, 2o)$.
		\item For all nodes $y\in\setY$ there exist three undirected edges $(y, y/2)$, $(y, 2y)$, $(y, o)$.
		\item For all vertices $v\in\setE-\setY$ there exist two undirected edges $(v, v/2)$, $(v, 2v)$. \Black
	\end{enumerate}
\end{theorem}

\section{Proof of the Collatz conjecture} 

\begin{figure}[htbp]
	\centering
	\includegraphics[scale=0.4]{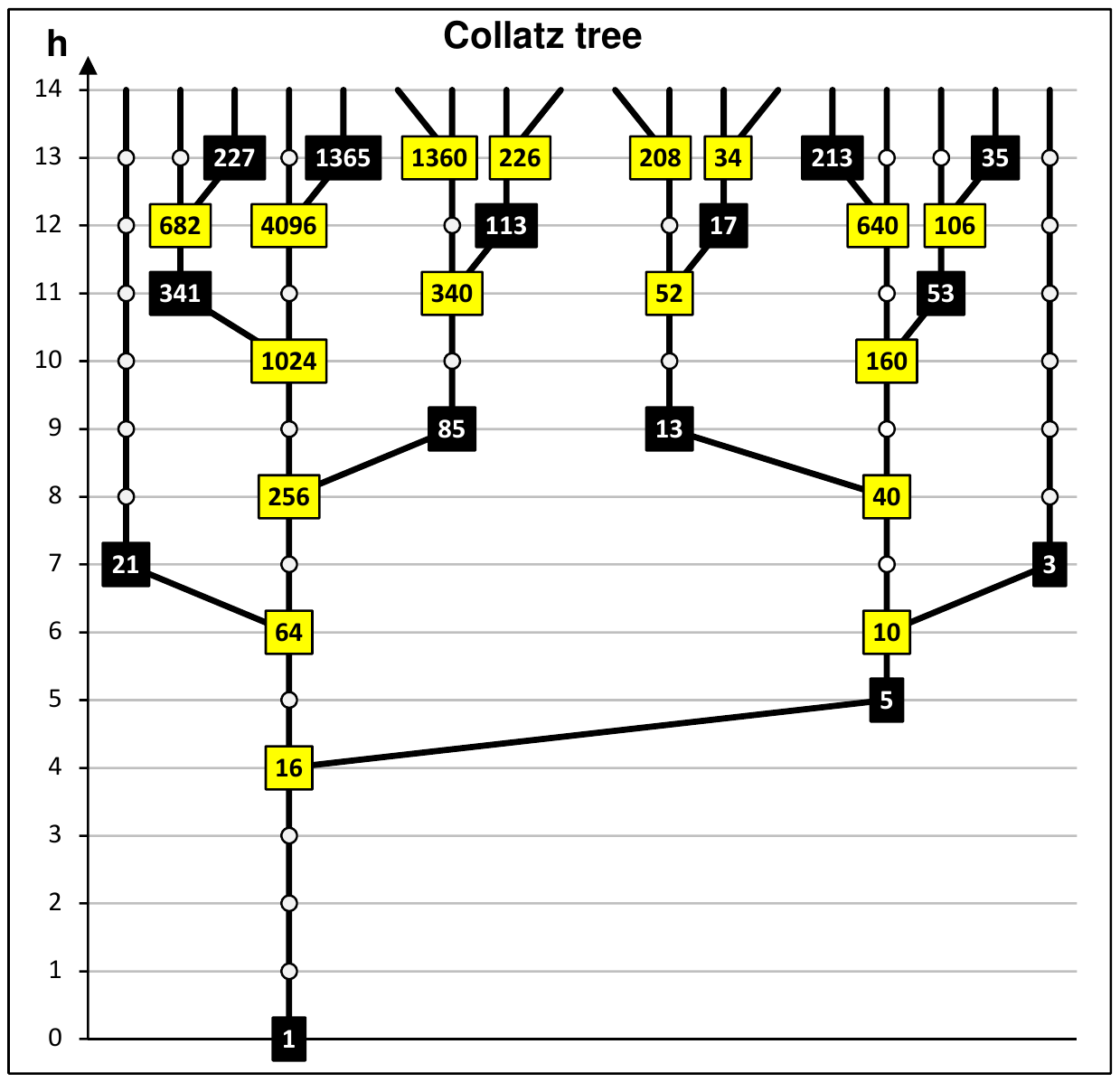}
	\caption{Height-oriented binary Collatz tree $\protect \Atr (1)$ up to level $h=13$ with an indicated continuation for $h=14$.}
	\label{fig:Ctv}
\end{figure}

\noindent The detour due to splitting the Collatz graph $\Cgr$ into separate components leads to important insights into the overall structure of this graph  provoked by the generating function. 
\begin{theorem}
The Collatz graph $\Cgr(1)$ is an infinite connected graph with vertex set $\setZ$.

\Proof 
\begin{enumerate}[nosep]
	\item \smallskip According to Theorems 5.1 and 5.3 is $V\bigl(\Cgr(1)\bigr)=\setZ$ guaranteeing the infiniteness too.
	\item \smallskip The inverse operations of $\Cor$ and $\Col$ provoke edges that match bijective maps. Thus $\Cgr(1)$ is an undirected graph.  
	\item \smallskip The graph $\Cgr(1)$ is connected. If we assume that it is not connected, there has to be at least one node $v\neq 1$ 
		which has no edge to a predecessor or successor. But this is a contradiction to the fact that the root $v=1$ is the only vertex 
		which has just one incident edge. All nodes $v\neq 1$ either have two or three definite incident edges according Item 2 and Theorem 5.4. \Black 
\end{enumerate}
\end{theorem}
\noindent The graph $\Cgr(1)$ is connected indeed but in the representation of (Figure~\ref{fig:Ct}) it seems to be an utter mess of edges crossing each other in an arbitrary manner. The existance of circuits cannot be excluded. On the contrary a binary tree is a well structured planar graph whose nodes can be arranged in height-oriented levels (Figure~\ref{fig:Ctv}).

\noindent If we are able to transform the graph $\Cgr(1)$ into a level-oriented binary tree 
		then this is a proof that there cannot be any circuits. 
		This could be done by a recursive procedure that evaluates the function $\Cor$ 
		and pulls all subtrees nearer to the tree they descent from. 
		But this fails because the recursion of any subtree will never come to an end 
		as there are no leaves that terminate the recursive descent. 
		However the function $\Cor$ offers an iterative procedure instead that assigns a height 
		and a dedicated level to all nodes of $\Cgr(1)$ {\footnotesize \fbox{2,3}}.
\begin{theorem}
The Collatz graph $\Cgr(1)$ can be transformed to an infinite undirected binary tree $\Atr(1)$.

\Proof We construct $\Cgr(1)$ by using the function $\Cor$ and induction. 
\begin{enumerate}[nosep]
	\item \smallskip We assume that $\Cgr(1)$ is a binary tree up to the level h=13 
		and as Figure~\ref{fig:Ctv} shows this is true for $h=13$.
	\item \smallskip In respect of $\Cor$ all nodes $v\ne1$ have one incoming edge only. 
		For even nodes $v$ these are the edges 
		$\overrightarrow e=\langle v/2, v\rangle$ and $\overrightarrow e=\langle y, o\rangle$ 
		for odd vertices $o\ne 1$ (Theorem 5.4).
	\item \smallskip $\Cor$ creates no nodes on level h that can have an outgoing edge 
		$\overrightarrow e=\langle v, 2v\rangle$ or $\overrightarrow e=\langle y,o\rangle$ 
		to a vertex on the levels from $0$ up to and including $h$ itself
		since all these nodes are already saturated regarding 
		to their indegree viz. the number of incoming edges.
	\item \smallskip Thus all successors of the nodes of the level $h$ have to be arranged 
		on the next higher level $h'=h+1$.
	\item \smallskip The constraints of the Items 2 to 4 apply to all nodes of every new level $h'$ and so 
		induction applies ad infinitum since $\Cgr(1)$ is connected.
\end{enumerate}

\smallskip \noindent There are no circuits in $\Cgr(1)$ thus it is an infinite binary tree $\Atr(1)$ 
		with vertex set $\setZ$ and therefore the Collatz conjecture is true. \Black
\end{theorem}

\section{References} 

\begin{enumerate}[nosep,label=\fbox{\arabic*}]
	\item Andaloro, Paul: The 3x+1 problem and directed graphs, Fibonacci Quarterly 40; 2002; p.43
	\item Andrei, S. et al.: Chains in Collatz's tree; Report 217; 1999; Department of Informatics; Universit\"{a}t Hamburg;
		 \url{http://edoc.sub.uni-hamburg.de/informatik/volltexte/2009/41/pdf/B\_217.pdf}
	\item Andrei, S. et al.: Some results on the Collatz problem; Acta Informatica 37; 2000; p.145
	\item Collatz, Lothar: Verzweigungsdiagramme und Hypergraphen; International Series for Numerical Mathematics; Vol.38; Birkkhäuser; 1977
	\item Diestel, R.: Graph Theory (GTM 137) $\mathrm{5}^{th}$edition; Springer-Verlag; New York; 2016
	\item Lang, W.: On Collatz' Words, Sequences and Trees; arXiv:1404.2710v1; 10 Apr 2014
	\item Laarhoven, Thijs \& de Weger, Benne: The Collatz conjecture and De Bruijn graphs; arXiv:1209.3495v1; 16 sep 2012
	\item Wirsching, G.: The Dynamical System Generated by the 3n+1 Function; Lecture Notes in Mathematics; Vol. 1681; Springer-Verlag; New York; 1998.
\end{enumerate}

\vspace{1cm}
\noindent Heinz Ebert, Im Heidgen 3, 53819 Neunkirchen-Seelscheid \hfill mail@heinzebert.de
\end{document}